\documentclass[11pt]{amsproc}

\usepackage{xcolor}
\usepackage{amsmath,amsthm,amssymb,mathtools,bm}
\allowdisplaybreaks
\usepackage{fullpage}
\usepackage{hyperref}
\usepackage[abbrev]{amsrefs}
\renewcommand{\MR}[1]{}
\usepackage{aliascnt}
\usepackage[capitalize,nameinlink,noabbrev,nosort]{cleveref}

\hypersetup{
  colorlinks=true,
  linkcolor=brown,
  citecolor=brown,
  filecolor=brown,
  urlcolor=brown,
  pdftitle={One of the numbers Z\_A(3) and Z\_A(5) is non-zero},
  pdfauthor={Toshiki Matsusaka}
}

\theoremstyle{plain}
\newaliascnt{theorem}{theoremcounter}
\newtheorem{theorem}[theorem]{Theorem}
\aliascntresetthe{theorem}
\newaliascnt{corollary}{theoremcounter}
\newtheorem{corollary}[corollary]{Corollary}
\aliascntresetthe{corollary}
\newaliascnt{proposition}{theoremcounter}
\newtheorem{proposition}[proposition]{Proposition}
\aliascntresetthe{proposition}
\newaliascnt{lemma}{theoremcounter}
\newtheorem{lemma}[lemma]{Lemma}
\aliascntresetthe{lemma}
\theoremstyle{definition}
\newaliascnt{definition}{theoremcounter}
\newtheorem{definition}[definition]{Definition}
\aliascntresetthe{definition}
\theoremstyle{remark}
\newaliascnt{remark}{theoremcounter}

\aliascntresetthe{remark}
\theoremstyle{plain}
\numberwithin{equation}{section}

\newcommand{\Z}{\mathbb{Z}}
\newcommand{\Q}{\mathbb{Q}}
\newcommand{\F}{\mathbb{F}}
\newcommand{\calA}{\mathcal{A}}
\newcommand{\calI}{\mathcal{I}}
\newcommand{\bk}{\bm{k}}
\newcommand{\bs}{\bm{s}}
\newcommand{\bu}{\bm{u}}
\newcommand{\bv}{\bm{v}}
\newcommand{\dd}{\mathrm{d}}
\DeclareMathOperator{\rank}{rank}
\DeclareMathOperator{\wt}{wt}
\DeclareMathOperator{\dep}{dep}

\begin{document}
\raggedbottom
\lefthyphenmin=2
\righthyphenmin=3

\title[Non-zero finite multiple zeta values]{One of the numbers $Z_\calA(3)$ and $Z_\calA(5)$ is non-zero}
\author[T.~Matsusaka]{Toshiki Matsusaka}
\address{Faculty of Mathematics, Kyushu University, Motooka 744, Nishi-ku, Fukuoka 819-0395, Japan}
\email{matsusaka@math.kyushu-u.ac.jp}

\begin{abstract}
We prove that one of the numbers $Z_\calA(3)$ and $Z_\calA(5)$ is non-zero, answering Kaneko and Zagier's question on the existence of non-zero finite multiple zeta values of positive weight.
For each integer $N\ge3$, we deduce that $\zeta_\calA(\bk)\ne0$ for some index $\bk$ with entries in $\{1,2\}$ and weight in $(N,2N]$.
\end{abstract}

\maketitle

\section{Introduction}\label{sec:introduction}

Finite multiple zeta values~\cite{KZ} are elements of the $\Q$-algebra
\[
 \calA\coloneqq\left(\prod_{p\text{ prime}}\Z/p\Z\right)
 \Big/\left(\bigoplus_{p\text{ prime}}\Z/p\Z\right).
\]
An index $\bk=(k_1,\ldots,k_r)$ is a tuple of positive integers. Its weight and depth are $\wt(\bk)\coloneqq k_1+\cdots+k_r$ and $\dep(\bk)\coloneqq r$, respectively. For an integer $N\ge1$, we put
\begin{align}\label{eq:truncated}
 \zeta_N(\bk)\coloneqq
 \sum_{0<n_1<\cdots<n_r<N}\frac{1}{n_1^{k_1}\cdots n_r^{k_r}},
 \qquad \zeta_N(\varnothing)=1,
\end{align}
where $\varnothing$ denotes the empty index, of weight and depth zero. By definition, the finite multiple zeta value $\zeta_\calA(\bk)\coloneqq(\zeta_p(\bk)\bmod p)_p\in\calA$ is zero if and only if $\zeta_p(\bk)\equiv0\pmod p$ for all sufficiently large primes. In depth one, we have $\zeta_p(k)\equiv0\pmod p$ for $p>k+1$, and hence $\zeta_\calA(k)=0$ for $k\ge1$. Kaneko and Zagier began their joint study of finite multiple zeta values in $\calA$ in 2012~\cite{KZ}*{Introduction, Remark~2}. Until now, it was not known whether the values of positive weight span a non-zero $\Q$-vector space~\cite{KZ}*{Introduction, Remark~1}.

Let $B_n$ denote the $n$-th Bernoulli number. Following~\cite{KZ}*{Section~1, Example~5}, we define
\begin{align}\label{eq:ZA}
 Z_\calA(k)\coloneqq\left(\frac{B_{p-k}}{k}\bmod p\right)_p
 \qquad(k\ge2),
\end{align}
where we omit the finitely many primes for which the component is undefined. By the double zeta formula~\cite{KZ}*{(48)}, we have $\zeta_\calA(1,2)=3Z_\calA(3)$ and $\zeta_\calA(1,4)=5Z_\calA(5)$.
\begin{theorem}\label{thm:main}
One of the numbers $Z_\calA(3)$ and $Z_\calA(5)$ is non-zero.
\end{theorem}

The theorem also answers Seki's Problem~1~\cite{Seki}*{Section~5} with $k\in\{3,5\}$. For a prime $p\ge5$, let $\eth_p$ denote the smallest odd integer $3\le k<p$ such that $B_{p-k}\not\equiv0\pmod p$. In this notation, the theorem is equivalent to $\eth_p\le5$ for infinitely many primes $p$.

Combining \cref{thm:main} with a generating function of Maesaka and Ono~\cite{MO}, we obtain infinitely many indices $\bk$ for which $\zeta_\calA(\bk)\ne0$.
\begin{corollary}\label{cor:large-weight}
For every integer $N\ge3$, there is an index $\bk$ with entries in $\{1,2\}$ such that $N<\wt(\bk)\le2N$ and $\zeta_\calA(\bk)\ne0$.
\end{corollary}

For $1\le\wt(\bk)\le6$, the value $\zeta_\calA(\bk)$ is a $\Q$-linear combination of $Z_\calA(3)$, $Z_\calA(5)$, and $Z_\calA(3)^2$~\cite{KZ}*{Section~1, Example~7}. Thus, if $Z_\calA(3)=Z_\calA(5)=0$, then every sufficiently large prime $p$ satisfies the vanishing condition
\begin{equation}\tag{\(\mathrm{V}_p\)}\label{eq:vanish}
 \zeta_p(\bk)\equiv0\pmod p
 \qquad\text{for every index }\bk\text{ with }1\le\wt(\bk)\le6.
\end{equation}

We introduce a matrix $K_m$ in \cref{sec:matrix} and, for $m=4\cdot17^r$, prove an unconditional upper bound for $|\det K_m|$ in \cref{sec:upper}. If~\eqref{eq:vanish} holds for all sufficiently large primes, the lower bound in \cref{sec:arithmetic} contradicts this upper bound.

\section{The matrix \texorpdfstring{$K_m$}{K\_m}}\label{sec:matrix}

For a prime $p$, we put $\Z_{(p)}=\{a/b:a,b\in\Z,\ p\nmid b\}$ and identify $\Z_{(p)}/p\Z_{(p)}$ with $\F_p$. Let $v_p$ denote the $p$-adic valuation with $v_p(p)=1$ and $v_p(0)=+\infty$. We write $[z^n]F(z)$ for the coefficient of $z^n$ in $F(z)$, and similarly for several variables.

For integers $j\ge0$, put $b_j(X)=\binom Xj\in\Q[X]$. The Hoffman dual $\bk^\vee$ of a non-empty index $\bk$ is obtained by writing its entries as sums of ones and interchanging plus signs and commas. We define the linear map $\Lambda_{\bk}:\Q[X]\to\Q$ on the binomial basis by
\begin{align}\label{eq:Lambda}
 \Lambda_{\bk}\bigl(b_j(X)\bigr)\coloneqq
 (-1)^{j+\wt(\bk)-1}\bigl(\zeta_{j+2}(\bk^\vee)-\zeta_{j+1}(\bk^\vee)\bigr),
 \qquad \Lambda_{\varnothing}=0.
\end{align}

\begin{lemma}\label{lem:moment}
Suppose that $p$ satisfies~\eqref{eq:vanish}. For $\wt(\bk)\le7$ and $P(X)\in\Z_{(p)}[X]$ of degree at most $p-2$, we have $\Lambda_{\bk}(P(X))\in\Z_{(p)}$ and
\begin{align}\label{eq:moment}
 \sum_{n=1}^{p}P(n-1)\zeta_n(\bk)=\Lambda_{\bk}\bigl(P(X)\bigr)
 \qquad\text{in }\F_p.
\end{align}
\end{lemma}

\begin{proof}
Since $j!$ is a unit in $\Z_{(p)}$ for $0\le j\le p-2$, the polynomials $b_0(X),\ldots,b_{p-2}(X)$ form a basis over $\Z_{(p)}$ in degrees at most $p-2$. By linearity, it suffices to take $P(X)=b_j(X)$.

By Pascal's identity, for $0\le a<p$ we have
\[
 \sum_{n=a+1}^{p}\binom{n-1}{j}
 =\sum_{n=a+1}^{p}\left(\binom n{j+1}-\binom{n-1}{j+1}\right)
 =\binom p{j+1}-\binom a{j+1}.
\]
For $\bk=\varnothing$, we have $\zeta_n(\varnothing)=1$, so taking $a=0$ gives $\binom p{j+1}$ for the left-hand side of~\eqref{eq:moment}. Since $1\le j+1<p$, this is $0=\Lambda_{\varnothing}(b_j(X))$ in $\F_p$.

Now let $\bk=(k_1,\ldots,k_r)$ be non-empty and put $w\coloneqq\wt(\bk)$. As each denominator in~\eqref{eq:Lambda} is a product of powers of integers at most $j+1<p$, we have $\Lambda_{\bk}(b_j(X))\in\Z_{(p)}$. For $a\ge1$, put $\Delta_a(\bk)\coloneqq\zeta_{a+1}(\bk)-\zeta_a(\bk)$.

Put $\overline{\bk}\coloneqq(k_r,\ldots,k_1)$. With the increasing order in~\eqref{eq:truncated}, our $\bk$ corresponds to $\overline{\bk}$ in Sakugawa--Seki~\cite{SS}*{Definition~1.2}. In our notation, the two polynomials in $\F_p[t]$ are
\[
\begin{aligned}
 \pounds_{p,\bk}(t)
 &\coloneqq\sum_{0<n_1<\cdots<n_r<p}
       \frac{t^{n_r}}{n_1^{k_1}\cdots n_r^{k_r}}
   =\sum_{a=1}^{p-1}\Delta_a(\bk)t^a,\\
 \widetilde{\pounds}^{\star}_{p,\bk}(t)
 &\coloneqq\sum_{0<n_1\le\cdots\le n_r<p}
       \frac{t^{n_1}}{n_1^{k_1}\cdots n_r^{k_r}}.
\end{aligned}
\]
Here $\pounds_{p,\bk}(1)=\zeta_p(\bk)$, and we write $\zeta_p^\star(\bk)=\widetilde{\pounds}^{\star}_{p,\bk}(1)$.

We follow the proof of~\cite{SS}*{Theorem~1.3} with $N=p-1$. First apply~\cite{SS}*{Theorem~2.11} with the index $\bk$ and $(t_1,\ldots,t_r)=(t,1,\ldots,1)$. For $1\le i<r$, each factor $\zeta_p^\star(k_{i+1},\ldots,k_r)$ is a sum of values $\zeta_p(\bs)$ with $1\le\wt(\bs)\le w-1\le6$. Hence the correction terms vanish by~\eqref{eq:vanish}, and
\[
 \widetilde{\pounds}^{\star}_{p,\bk}(t)=(-1)^{r-1}\pounds_{p,\overline{\bk}}(t).
\]
Next apply~\cite{SS}*{Corollary~2.7} with the index $\bk$ and the variable $1-t$. Since $(-1)^n\binom{p-1}{n}=1$ in $\F_p$ for $0\le n<p$, we obtain
\begin{align*}
 \pounds_{p,\bk}(1-t)
 &=(-1)^{r-1}\widetilde{\pounds}^{\star}_{p,\overline{\bk}}(1-t)\\
 &=(-1)^{r-1}\bigl(\widetilde{\pounds}^{\star}_{p,\overline{\bk^\vee}}(t)
                        -\zeta_p^\star(\overline{\bk^\vee})\bigr)\\
 &=(-1)^{w-1}\bigl(\pounds_{p,\bk^\vee}(t)-\zeta_p(\bk^\vee)\bigr).
\end{align*}
Here $\dep(\bk)+\dep(\bk^\vee)=\wt(\bk)+1$.

We now compute the left-hand side of~\eqref{eq:moment} for $P(X)=b_j(X)$. Using $[t^{j+1}]\zeta_p(\bk^\vee)=0$ for $j\ge0$, we obtain in $\F_p$
\begin{align*}
 \sum_{n=1}^p\binom{n-1}j\zeta_n(\bk)
 &=\sum_{a=1}^{p-1}\Delta_a(\bk)\sum_{n=a+1}^{p}\binom{n-1}j\\
 &=-\sum_{a=1}^{p-1}\binom a{j+1}\Delta_a(\bk)
   =(-1)^j[t^{j+1}]\pounds_{p,\bk}(1-t)\\
 &=(-1)^{j+w-1}[t^{j+1}]\pounds_{p,\bk^\vee}(t)\\
 &=(-1)^{j+w-1}\bigl(\zeta_{j+2}(\bk^\vee)-\zeta_{j+1}(\bk^\vee)\bigr)
   =\Lambda_{\bk}\bigl(b_j(X)\bigr).\qedhere
\end{align*}
\end{proof}

For $\bk=(k_1,\ldots,k_d)\ne\varnothing$ and $|z|<1$, the multiple polylogarithm is defined by
\[
 \operatorname{Li}_{\bk}(z)\coloneqq
 \sum_{0<n_1<\cdots<n_d}\frac{z^{n_d}}{n_1^{k_1}\cdots n_d^{k_d}}
 =\sum_{n\ge1}\bigl(\zeta_{n+1}(\bk)-\zeta_n(\bk)\bigr)z^n.
\]
We put $\operatorname{Li}_{\varnothing}(z)\coloneqq1$. With $z=1-e^t$, termwise differentiation gives
\[
\begin{aligned}
 \frac{\dd}{\dd t}\operatorname{Li}_{(\bs,1)}(1-e^t)
 &=-\operatorname{Li}_{\bs}(1-e^t),\\
 \frac{\dd}{\dd t}\operatorname{Li}_{(\bs,a+1)}(1-e^t)
 &=\frac{e^t}{e^t-1}\operatorname{Li}_{(\bs,a)}(1-e^t)\qquad(a\ge1).
\end{aligned}
\]
The first equality uses $\sum_{n>h}z^{n-1}=z^h/(1-z)$. Applying $\Lambda_{\bk}$ to the coefficients of $e^{tX}=\sum_{j\ge0}b_j(X)(e^t-1)^j$ and using~\eqref{eq:Lambda}, we obtain near $t=0$
\begin{equation}\label{eq:polylog}
 F_{\bk}(t)\coloneqq\sum_{n\ge0}\Lambda_{\bk}(X^n)\frac{t^n}{n!}
 =(-1)^{\wt(\bk)-1}\frac{\operatorname{Li}_{\bk^\vee}(1-e^t)}{1-e^t},
 \qquad F_{\varnothing}=0.
\end{equation}

We first describe the row and column indices of the matrices $G_{p,m}$ and $K_m$ in \cref{def:K}. For an integer $d\ge0$, let $\calI_d$ denote the set of indices of weight at most $d$, including $\varnothing$. For each integer $m\ge1$, we take
\[
 \text{rows }(i,\bu)\in\{0,\ldots,2m-1\}\times\calI_3,\qquad
 \text{columns }(j,\bv)\in\{0,\ldots,m-1\}\times\calI_4,
\]
ordered lexicographically by $(i,\wt(\bu),\bu)$ and $(j,\wt(\bv),\bv)$, respectively. Since $|\calI_3|=8$ and $|\calI_4|=16$, both matrices have size $16m\times16m$.

Let $*$ denote the harmonic (stuffle) product on formal $\Q$-linear combinations of indices, with unit $\varnothing$~\cite{KZ}*{Section~7(a)}. We extend $\zeta_N$ and $\bk\mapsto\Lambda_{\bk}$ by linearity. By the definition of $*$, we have
\begin{align}\label{eq:shuffle}
 \zeta_N(\bu*\bv)=\zeta_N(\bu)\zeta_N(\bv).
\end{align}

\begin{definition}\label{def:K}
For a prime $p$, we define $G_{p,m}\in M_{16m}(\F_p)$ and $K_m\in M_{16m}(\Q)$ by
\begin{align}
 (G_{p,m})_{(i,\bu),(j,\bv)}
 &\coloneqq\sum_{n=1}^{p}b_i(n-1)b_j(n-1)\zeta_n(\bu)\zeta_n(\bv),\label{eq:G}\\*
 (K_m)_{(i,\bu),(j,\bv)}
 &\coloneqq\Lambda_{\bu*\bv}\bigl(b_i(X)b_j(X)\bigr).\label{eq:matrix}
\end{align}
\end{definition}

\begin{proposition}\label{prop:reduction}
We have $\rank_{\F_p}G_{p,m}\le p$. If $p$ satisfies~\eqref{eq:vanish} and $i+j\le p-2$, then $(K_m)_{(i,\bu),(j,\bv)}\in\Z_{(p)}$ and
\[
 (K_m)_{(i,\bu),(j,\bv)}=(G_{p,m})_{(i,\bu),(j,\bv)}
 \qquad\text{in }\F_p.
\]
\end{proposition}

\begin{proof}
Over $\F_p$, we have
\[
 G_{p,m}
 =\underbrace{\bigl(b_i(n-1)\zeta_n(\bu)\bigr)_{(i,\bu),n}}_{16m\times p}
  \underbrace{\bigl(b_j(n-1)\zeta_n(\bv)\bigr)_{n,(j,\bv)}}_{p\times16m}
 \qquad(1\le n\le p)
\]
and hence $\rank_{\F_p}G_{p,m}\le p$. Suppose that $p$ satisfies~\eqref{eq:vanish} and $i+j\le p-2$. Then $b_i(X)b_j(X)\in\Z_{(p)}[X]$ has degree at most $p-2$, and every index in $\bu*\bv$ has weight at most seven. By \cref{lem:moment}, the corresponding entry of $K_m$ is $p$-integral and, using~\eqref{eq:shuffle}, we obtain
\[
 (K_m)_{(i,\bu),(j,\bv)}
 =\sum_{n=1}^{p}b_i(n-1)b_j(n-1)\zeta_n(\bu*\bv)
 =(G_{p,m})_{(i,\bu),(j,\bv)}
 \qquad\text{in }\F_p.\qedhere
\]
\end{proof}

\begin{samepage}
By the binomial theorem,
\[
 (1+x)^X(1+y)^X=(1+x+y+xy)^X
 =\sum_{l\ge0}b_l(X)(x+y+xy)^l.
\]
Comparing the coefficients of $x^iy^j$ and using~\eqref{eq:matrix}, we obtain
\begin{subequations}
\begin{gather}
 b_i(X)b_j(X)=\sum_{l=\max(i,j)}^{i+j}C_{ij}^{l}b_l(X),
 \qquad
 C_{ij}^{l}=[x^iy^j](x+y+xy)^l=\binom li\binom i{i+j-l},\label{eq:binprod}\\
 (K_m)_{(i,\bu),(j,\bv)}
 =\Lambda_{\bu*\bv}\bigl(b_i(X)b_j(X)\bigr)
 =\sum_{l=\max(i,j)}^{i+j}C_{ij}^{l}\Lambda_{\bu*\bv}\bigl(b_l(X)\bigr).\label{eq:matrix-expansion}
\end{gather}
\end{subequations}
The formula for $C_{ij}^l$ follows by writing $x+y+xy=y+x(1+y)$.
\par\end{samepage}

\section{An upper bound for \texorpdfstring{$|\det K_m|$}{|det K\_m|}}\label{sec:upper}

\begin{lemma}\label{lem:scaling}
We have $K_4\in M_{64}(\Z_{(17)})$ and $\det K_4\not\equiv0\pmod{17}$. Let $q=17^r$ with $r\in\Z_{\ge0}$. Multiply row $(i,\bu)$ of $K_{4q}$ by $q^{\wt(\bu)}$ and column $(j,\bv)$ by $q^{\wt(\bv)}$. The resulting matrix has entries in $\Z_{(17)}$ and is non-singular modulo $17$. More precisely, we have
\[
 q^{\wt(\bu)+\wt(\bv)}(K_{4q})_{(aq+s,\bu),(bq+t,\bv)}
 =(K_4)_{(a,\bu),(b,\bv)}C_{st}^{q-1}
 \qquad\text{in }\F_{17}.
\]
Here $0\le a<8$, $0\le b<4$, and $0\le s,t<q$.
\end{lemma}

\begin{proof}
In~\eqref{eq:matrix-expansion}, take $m=4q$, $i=aq+s$, and $j=bq+t$. Fix $\max(i,j)\le l\le i+j$.
For a non-empty index $\bk$ with $\bk^\vee=(\gamma_1,\ldots,\gamma_d)$, it follows from~\eqref{eq:Lambda} that
\[
 q^{\wt(\bk)}\Lambda_{\bk}\bigl(b_l(X)\bigr)
 =(-1)^{l+\wt(\bk)-1}\sum_{0<n_1<\cdots<n_d=l+1}
   \prod_{\alpha=1}^{d}\left(\frac{q}{n_\alpha}\right)^{\gamma_\alpha}.
\]
Since $n_\alpha\le l+1\le i+j+1<12q<17q$, each $q/n_\alpha$ is $17$-integral. Every index in $\bu*\bv$ has weight $\wt(\bu)+\wt(\bv)$, so the entries after multiplication are integer linear combinations of these sums. A direct calculation gives $\det K_4\equiv6\pmod{17}$.

Only terms with $q\mid n_\alpha$ for every $\alpha$ contribute modulo $17$. If $l+1=q(c+1)$, dividing each $n_\alpha$ by $q$ gives the sum for $\Lambda_{\bk}(b_c(X))$, with the same sign since $l-c=(q-1)(c+1)$ is even. Hence
\begin{equation}\label{eq:scaling}
 q^{\wt(\bk)}\Lambda_{\bk}\bigl(b_l(X)\bigr)
 \equiv
 \begin{cases}
  \Lambda_{\bk}\bigl(b_c(X)\bigr) & \text{if } l+1=q(c+1),\\
  0 & \text{if } q\nmid l+1
 \end{cases}
 \pmod{17}.
\end{equation}
For the empty index, both sides are zero.

For $c\ge0$, the identity $(x+y+xy)^q=x^q+y^q+x^qy^q$ in $\F_{17}[x,y]$ gives
\[
 (x+y+xy)^{cq+q-1}=(x^q+y^q+x^qy^q)^c(x+y+xy)^{q-1}.
\]
The last factor has degree less than $q$ in each variable, so it can contribute only $x^sy^t$ to the coefficient of $x^{aq+s}y^{bq+t}$. Thus $C_{ij}^{cq+q-1}=C_{ab}^c C_{st}^{q-1}$. Since $C_{ab}^c=[x^ay^b](x+y+xy)^c$ vanishes unless $\max(a,b)\le c\le a+b$, the scaled entry in row $(i,\bu)$ and column $(j,\bv)$ of $K_{4q}$ is, by~\eqref{eq:matrix-expansion} and~\eqref{eq:scaling},
\[
 C_{st}^{q-1}\sum_{c=\max(a,b)}^{a+b}C_{ab}^c\Lambda_{\bu*\bv}\bigl(b_c(X)\bigr)
 =C_{st}^{q-1}(K_4)_{(a,\bu),(b,\bv)}
 \qquad\text{in }\F_{17}.
\]
Reordering the rows by $((a,\bu),s)$ and the columns by $((b,\bv),t)$ gives $K_4\otimes(C_{st}^{q-1})_{0\le s,t<q}$ over $\F_{17}$. The second factor is anti-triangular since $C_{st}^{q-1}=0$ for $s+t<q-1$. Its anti-diagonal entries are $C_{s,q-1-s}^{q-1}=\binom{q-1}s=(-1)^s$ in $\F_{17}$. Thus its determinant is $(-1)^{\binom q2}\prod_{s=0}^{q-1}(-1)^s=1$ in $\F_{17}$, and both factors are non-singular.
\end{proof}

\begin{proposition}\label{prop:matrix}
For $m=4\cdot17^r$, $r\in\Z_{\ge0}$, we have $\det K_m\ne0$ and, as $r\to\infty$,
\begin{align}\label{eq:analytic}
 \log|\det K_m|\le\bigl(-24\log4+o(1)\bigr)m^2.
\end{align}
\end{proposition}

\begin{proof}
By \cref{lem:scaling}, $\det K_m\ne0$. In the strip $|\operatorname{Im}t|<2\pi$, the only zero of $1-e^t$ is simple and occurs at $0$, where $\operatorname{Li}_{\bk}(1-e^t)$ vanishes for every non-empty $\bk$. The singularities at $0$ in the derivative formulas preceding~\eqref{eq:polylog} are therefore removable. Integrating from $0$ in this simply connected strip extends $\operatorname{Li}_{\bk}(1-e^t)$ holomorphically by induction on the weight, and~\eqref{eq:polylog} gives the same for $F_{\bk}$.

For $|z|<4$, put
\[
 g(z)\coloneqq2\log\frac{1+z/4}{1-z/4},\quad g(0)=0,\qquad
 Q_i(X)\coloneqq[z^i]e^{Xg(z)}.
\]
The quotient $(1+z/4)/(1-z/4)$ has positive real part, so $|\operatorname{Im}g(z)|<\pi$. Since $g(z)=z+O(z^2)$, $Q_i(X)-b_i(X)$ has degree less than $i$. Replacing $b_i(X)$ by $Q_i(X)$ adds a combination of earlier rows with the same $\bu$, and likewise for columns with the same $\bv$. This preserves the determinant. By~\eqref{eq:polylog}, the new entries are
\[
 \Lambda_{\bu*\bv}\bigl(Q_i(X)Q_j(X)\bigr)
 =[z^iw^j]F_{\bu*\bv}\bigl(g(z)+g(w)\bigr).
\]
For $0<R<4$, Cauchy's estimate on $|z|=|w|=R$ gives $O_R(R^{-i-j})$ since $|\operatorname{Im}(g(z)+g(w))|<2\pi$. Multiplying row $(i,\bu)$ by $R^i$ and column $(j,\bv)$ by $R^j$ gives row norms at most $C_R\sqrt m$. The total exponent of $R$ in the determinant is $8\binom{2m}{2}+16\binom m2=24m^2-16m$. Hadamard's inequality gives
\[
 \frac{\log|\det K_m|}{m^2}\le-\left(24-\frac{16}{m}\right)\log R
       +\frac{16}{m}\log(C_R\sqrt m).
\]
For fixed $0<R<4$ the last term tends to zero as $r\to\infty$. Taking the upper limit and then letting $R\uparrow4$ gives~\eqref{eq:analytic}.
\end{proof}

\section{A lower bound for \texorpdfstring{$|\det K_m|$}{|det K\_m|}}\label{sec:arithmetic}

We first bound the denominators of the entries of $K_m$ at a prime $p$. For an index $\bs$, let $\nu(\bs)$ denote the number of entries of $\bs$ equal to one. For real $y$, put $y_+=\max(y,0)$. We write $\mathbf1_A=1$ if the condition $A$ holds and $\mathbf1_A=0$ otherwise.

\begin{lemma}\label{lem:denominator}
Suppose that $p^2>3m$ and that $p$ satisfies~\eqref{eq:vanish}. For each row $(i,\bu)$ and column $(j,\bv)$ of $K_m$, put
\begin{equation}\label{eq:local-exponents}
\begin{split}
 r_{i,\bu}&=\nu(\bu)\mathbf1_{i+m\ge p}+\lfloor i/p\rfloor
       +\bigl(\lfloor(i+m)/p\rfloor-1\bigr)_+,\\
 t_{j,\bv}&=\nu(\bv)\mathbf1_{j+2m\ge p}+\lfloor j/p\rfloor.
\end{split}
\end{equation}
Then $v_p((K_m)_{(i,\bu),(j,\bv)})\ge-r_{i,\bu}-t_{j,\bv}$.
\end{lemma}

\begin{proof}
Fix a row $(i,\bu)$ and a column $(j,\bv)$ in~\eqref{eq:matrix-expansion}. For $\max(i,j)\le l\le i+j$ and a non-empty index $\bs$ in $\bu*\bv$, we estimate $v_p\bigl(C_{ij}^l\Lambda_{\bs}(b_l(X))\bigr)$. Write $\bs^\vee=(\gamma_1,\ldots,\gamma_d)$. For each tuple $0<n_1<\cdots<n_d=l+1$, put $J=\{1\le a\le d : p\mid n_a\}$. By~\eqref{eq:Lambda} applied to $b_l(X)$,
\[
\begin{gathered}
 (-1)^{l+\wt(\bs)-1}\Lambda_{\bs}\bigl(b_l(X)\bigr)
 =\zeta_{l+2}(\bs^\vee)-\zeta_{l+1}(\bs^\vee)
 =\sum_{J\subseteq\{1,\ldots,d\}}H_J,\\
 H_J\coloneqq
 \sum_{\substack{0<n_1<\cdots<n_d=l+1\\
        p\mid n_a\ (a\in J),\ p\nmid n_a\ (a\notin J)}}
 \frac{1}{n_1^{\gamma_1}\cdots n_d^{\gamma_d}}.
\end{gathered}
\]
Since $i\le2m-1$ and $j\le m-1$, we have $n_a\le l+1\le i+j+1\le3m-1<p^2$. Thus the denominator $n_1^{\gamma_1}\cdots n_d^{\gamma_d}$ of every summand of $H_J$ has valuation $\sum_{a\in J}\gamma_a$. In particular, $H_{\varnothing}\in\Z_{(p)}$.

Suppose that $J\ne\varnothing$ and the sum defining $H_J$ is non-empty. Put $h=|J|$ and $b=\min J$. We claim that
\[
 v_p(H_J)\ge-\sum_{a\in J}\gamma_a+\mathbf1_{1\notin J}.
\]
For $b=1$ this follows from the valuation of each denominator. For $b>1$, take any $n_b,\ldots,n_d$ occurring in $H_J$ and write $n_b=cp$. The sum over $n_1,\ldots,n_{b-1}$ satisfies
\[
 \sum_{\substack{0<n_1<\cdots<n_{b-1}<cp\\p\nmid n_1\cdots n_{b-1}}}
 \frac{1}{n_1^{\gamma_1}\cdots n_{b-1}^{\gamma_{b-1}}}
 \equiv\sum_{0=d_0\le\cdots\le d_c=b-1}
 \prod_{t=1}^{c}\zeta_p(\gamma_{d_{t-1}+1},\ldots,\gamma_{d_t})\equiv0\pmod p.
\]
Here $d_t\coloneqq\#\{1\le a<b:n_a<tp\}$. Subtracting $(t-1)p$ from the $n_a$ in $((t-1)p,tp)$ gives the $t$-th factor, with $\zeta_p(\varnothing)=1$ if $d_t=d_{t-1}$. Since $b-1>0$, each product has a factor of positive weight at most $\gamma_1+\cdots+\gamma_{b-1}\le\wt(\bs)-1\le6$, which vanishes by~\eqref{eq:vanish}. The inner sum thus contributes a factor of $p$, while the denominator $n_b^{\gamma_b}\cdots n_d^{\gamma_d}$ has valuation $\sum_{a\in J}\gamma_a$.

The exponent $\sum_{a\in J}\gamma_a$ satisfies
\begin{equation}\label{eq:dual-ones}
\begin{split}
 \sum_{a\in J}\gamma_a
 &\le\nu(\bs)+\sum_{a\in J}(\mathbf1_{a>1}+\mathbf1_{a<d})
  =\nu(\bs)+2h-\mathbf1_{1\in J}-\mathbf1_{d\in J}\\
 &\le\nu(\bs)+h+\lfloor l/p\rfloor-\mathbf1_{1\in J}.
\end{split}
\end{equation}
To prove the first inequality in~\eqref{eq:dual-ones}, write each $\gamma_a$ as a sum of ones. Taking the Hoffman dual interchanges plus signs and commas. For example,
\[
 (2,1,2)=(1+1,\,1,\,1+1)
 \quad\longmapsto\quad
 (1,\,1+1+1,\,1)=(1,3,1).
\]
If $a>1$, the first $1$ in $\gamma_a$ is added to the $1$ on its left. If $a<d$, the last $1$ is added to the $1$ on its right. Thus $\gamma_a$ gives $(\gamma_a-\mathbf1_{a>1}-\mathbf1_{a<d})_+$ entries equal to $1$ in $\bs$. The total over $a\in J$ is at most $\nu(\bs)$.

The second inequality follows from $h-\mathbf1_{d\in J}=\#\{1\le a<d : p\mid n_a\}\le\lfloor l/p\rfloor$, since $n_1,\ldots,n_{d-1}$ are distinct integers in $\{1,\ldots,l\}$.

For the coefficient $C_{ij}^l$, we have
\[
 v_p(C_{ij}^l)\ge v_p\!\left(\binom li\right)
 =\lfloor l/p\rfloor-\lfloor i/p\rfloor-\lfloor(l-i)/p\rfloor
 \ge\lfloor l/p\rfloor-\lfloor i/p\rfloor-\lfloor j/p\rfloor .
\]
Indeed, $\binom li$ divides $C_{ij}^l$ by~\eqref{eq:binprod}, and $v_p(n!)=\lfloor n/p\rfloor$ for $0\le n\le l<p^2$. Also, $l-i\le j$.

Combining the bounds for $v_p(H_J)$ and $v_p(C_{ij}^l)$ with~\eqref{eq:dual-ones}, we obtain
\[
 \begin{aligned}
 v_p(C_{ij}^l H_J)
 &\ge v_p(C_{ij}^l)-\sum_{a\in J}\gamma_a+\mathbf1_{1\notin J}\\
 &\ge 1-\nu(\bs)-h-\lfloor i/p\rfloor-\lfloor j/p\rfloor\\
 &\ge-r_{i,\bu}-t_{j,\bv}.
 \end{aligned}
\]
Here the terms $\lfloor l/p\rfloor$ cancel and $\mathbf1_{1\in J}+\mathbf1_{1\notin J}=1$. Since $n_1,\ldots,n_d$ include $h$ distinct multiples of $p$, we have $p\le hp\le l+1\le i+j+1\le\min(i+m,j+2m)$. Both indicators in~\eqref{eq:local-exponents} therefore equal one, and $h-1\le(\lfloor(i+m)/p\rfloor-1)_+$. Merging two entries in the stuffle product gives an entry at least two, so $\nu(\bs)\le\nu(\bu)+\nu(\bv)$. Since $\Lambda_{\varnothing}=0$, the entry $(K_m)_{(i,\bu),(j,\bv)}$ is an integer linear combination of the $C_{ij}^lH_J$. This proves the lemma.
\end{proof}

\begin{lemma}\label{lem:weighted-prime-sum}
For fixed $0<a<b$, as $m\to\infty$ we have
\[
 \frac1{m^2}\sum_{p>am}(bm-p)_+\log p=\frac{(b-a)^2}{2}+o(1).
\]
\end{lemma}

\begin{proof}
Put $\vartheta(y)\coloneqq\sum_{p\le y}\log p$. For $am<p\le bm$, write $bm-p=\int_p^{bm}\dd t$. Interchanging the sum and the integral and then putting $t=mx$ gives
\[
 \frac1{m^2}\sum_{p>am}(bm-p)_+\log p
 =\frac1{m^2}\int_{am}^{bm}\sum_{am<p\le t}\log p\,\dd t
 =\int_a^b\frac{\vartheta(mx)-\vartheta(am)}{m}\,\dd x.
\]
Since $mx\ge am\to\infty$ and $x\le b$, the prime number theorem gives $\vartheta(mx)/m=x\,\vartheta(mx)/(mx)=x+o(1)$ uniformly for $a\le x\le b$. The integral therefore tends to $\int_a^b(x-a)\,\dd x=(b-a)^2/2$.
\end{proof}

For the rest of this section, let $m=4\cdot17^r$ with $r\in\Z_{\ge0}$. For $p>3m/4$ we have $i,j,i+m<3m<4p$, so each floor in~\eqref{eq:local-exponents} is at most three. We therefore scale the rows and columns by setting
\[
 (K_m^{(p)})_{(i,\bu),(j,\bv)}\coloneqq
 \begin{cases}
 p^{h\wt(\bu)}(K_m)_{(i,\bu),(j,\bv)}p^{h\wt(\bv)}
    & \text{if }p\le3m/4,\quad h\coloneqq\lfloor\log_p(3m)\rfloor,\\
 p^{r_{i,\bu}}(K_m)_{(i,\bu),(j,\bv)}p^{t_{j,\bv}}
    & \text{if }p>3m/4,\quad r_{i,\bu},\,t_{j,\bv}\text{ as in~\eqref{eq:local-exponents}}.
 \end{cases}
\]
Let $E_p$ be the sum of the row and column exponents, so $\det K_m^{(p)}=p^{E_p}\det K_m$. For $p\ge3m$ all exponents vanish. Since $\det K_m\in\Q^\times$ by \cref{prop:matrix}, both $\sum_p E_p\log p$ and $\sum_p v_p(\det K_m^{(p)})\log p$ have only finitely many non-zero terms.

\begin{lemma}\label{lem:total-exponents}
Suppose that~\eqref{eq:vanish} holds for all sufficiently large primes. For all sufficiently large $r$, the matrix $K_m^{(p)}$ has entries in $\Z_{(p)}$ for every prime $p$. Moreover, as $r\to\infty$,
\[
 \sum_p E_p\log p=\left(\frac{519}{4}+o(1)\right)m^2.
\]
\end{lemma}

\begin{proof}
Take $r$ sufficiently large. Suppose first that $p\le3m/4$. By~\eqref{eq:matrix-expansion} and~\eqref{eq:Lambda}, the denominators in $(K_m)_{(i,\bu),(j,\bv)}$ have the form $n_1^{\gamma_1}\cdots n_d^{\gamma_d}$, where $\bs^\vee=(\gamma_1,\ldots,\gamma_d)$ for a non-empty index $\bs$ in $\bu*\bv$. Since $n_a\le n_d=l+1\le i+j+1<3m<p^{h+1}$, we have $v_p(n_a)\le h$. Thus
\[
 v_p\bigl(n_1^{\gamma_1}\cdots n_d^{\gamma_d}\bigr)
 =\sum_{a=1}^d\gamma_a v_p(n_a)
 \le h\sum_{a=1}^d\gamma_a
 =h\bigl(\wt(\bu)+\wt(\bv)\bigr).
\]
Here $\sum_a\gamma_a=\wt(\bs^\vee)=\wt(\bu)+\wt(\bv)$. Multiplication by $p^{h(\wt(\bu)+\wt(\bv))}$ therefore makes each summand $p$-integral.

There are $2m$ rows for each $\bu\in\calI_3$ and $m$ columns for each $\bv\in\calI_4$. Thus
\[
 E_p=2mh\sum_{\bu\in\calI_3}\wt(\bu)+mh\sum_{\bv\in\calI_4}\wt(\bv)=83mh.
\]

We have $h\log p\le\log(3m)$ since $p^h\le3m$, and $h=1$ for $\sqrt{3m}<p\le3m/4$ since $p\le3m<p^2$. Hence $\sum_{p\le3m/4}(h-1)\log p\le\sqrt{3m}\log(3m)$. The prime number theorem gives
\[
 \frac1{m^2}\sum_{p\le3m/4}E_p\log p
 =\frac{83}{m}\sum_{p\le3m/4}\log p+O\!\left(\frac{\log m}{\sqrt m}\right)
 =\frac{249}{4}+o(1).
\]

Now let $p>3m/4$. For large $m$, we have $p^2>3m$ and~\eqref{eq:vanish} holds, so \cref{lem:denominator} gives $K_m^{(p)}\in M_{16m}(\Z_{(p)})$. In~\eqref{eq:local-exponents}, $i<2m<3p$, $j<m<2p$, and $i+m<3m<4p$. For example, $(\lfloor(i+m)/p\rfloor-1)_+=\mathbf1_{i+m\ge2p}+\mathbf1_{i+m\ge3p}$. Expanding the floors gives
\[
\begin{aligned}
 \sum_{\bu\in\calI_3}r_{i,\bu}
 &=8\bigl(\mathbf1_{i+m\ge p}+\mathbf1_{i\ge p}+\mathbf1_{i\ge2p}
             +\mathbf1_{i+m\ge2p}+\mathbf1_{i+m\ge3p}\bigr),\\
 \sum_{\bv\in\calI_4}t_{j,\bv}
 &=20\mathbf1_{j+2m\ge p}+16\mathbf1_{j\ge p}.
\end{aligned}
\]
Here $|\calI_3|=8$, $|\calI_4|=16$, $\sum_{\bu\in\calI_3}\nu(\bu)=8$, and $\sum_{\bv\in\calI_4}\nu(\bv)=20$. For integers $A\le B$ and $c$, we have $\sum_{k=A}^{B-1}\mathbf1_{k\ge c}=(B-c)_+-(A-c)_+$. For the first term in the row sum, this gives
\[
 \sum_{i=0}^{2m-1}\mathbf1_{i+m\ge p}
 =\sum_{k=m}^{3m-1}\mathbf1_{k\ge p}=(3m-p)_+-(m-p)_+.
\]
Summing the row exponents over $0\le i<2m$ and the column exponents over $0\le j<m$, we obtain
\begin{equation}\label{eq:exponent-total}
\begin{split}
 E_p
 &=8\bigl((3m-p)_+-(m-p)_++(2m-p)_++(2m-2p)_++(3m-2p)_++(3m-3p)_+\bigr)\\
 &\quad+20\bigl((3m-p)_+-(2m-p)_+\bigr)+16(m-p)_+\\
 &=28(3m-p)_+-12(2m-p)_++48(m-p)_++8(3m-2p)_+.
\end{split}
\end{equation}
Since $(3m-2p)_+=2(3m/2-p)_+$, applying \cref{lem:weighted-prime-sum} to the four terms of~\eqref{eq:exponent-total} with $a=3/4$ and $b=3,2,1,3/2$, respectively, gives
\[
 \lim_{r\to\infty}\dfrac{1}{m^2}\sum_{p>3m/4}E_p\log p
 =\dfrac{1}{2}\left[28\left(\dfrac{9}{4}\right)^2-12\left(\dfrac{5}{4}\right)^2
                 +48\left(\dfrac{1}{4}\right)^2+16\left(\dfrac{3}{4}\right)^2\right]
 =\dfrac{135}{2}.
\]
Adding the contributions from the two ranges, we obtain
\begin{equation}\label{eq:denominator-total}
 \frac1{m^2}\sum_p E_p\log p
 =\frac{249}{4}+\frac{135}{2}+o(1)=\frac{519}{4}+o(1).\qedhere
\end{equation}
\end{proof}

\begin{lemma}\label{lem:scaled-divisibility}
Suppose that~\eqref{eq:vanish} holds for all sufficiently large primes. As $r\to\infty$,
\[
 \sum_p v_p(\det K_m^{(p)})\log p\ge\left(\frac{193}{2}+o(1)\right)m^2.
\]
\end{lemma}

\begin{proof}
For sufficiently large $r$, the entries of $K_m^{(p)}$ lie in $\Z_{(p)}$ for every prime $p$ by \cref{lem:total-exponents}. Put $\rho\coloneqq\rank_{\F_p}(K_m^{(p)}\bmod p)$ and choose rows $R_1,\ldots,R_\rho$ whose reductions form a basis of the row space modulo $p$. For each remaining row $R$, choose integers $c_j$ representing the coefficients of $R\bmod p$ in this basis. Then
\[
 R\equiv\sum_{j=1}^{\rho}c_jR_j\pmod p,
 \qquad R-\sum_{j=1}^{\rho}c_jR_j\in\bigl(p\Z_{(p)}\bigr)^{16m}.
\]
Subtracting $\sum_j c_jR_j$ from $R$ preserves the determinant. Factoring $p$ from each of the $16m-\rho$ resulting rows gives $v_p(\det K_m^{(p)})\ge16m-\rho$.

Now let $p>m$, so that~\eqref{eq:vanish} holds. If $i<p-m$, then~\eqref{eq:local-exponents} gives $r_{i,\bu}=0$, and $i+j\le p-2$ for every $j<m$. By \cref{prop:reduction},
\[
 (K_m^{(p)})_{(i,\bu),(j,\bv)}
 =p^{t_{j,\bv}}(G_{p,m})_{(i,\bu),(j,\bv)}
 \qquad\text{in }\F_p.
\]
The rows of $K_m^{(p)}\bmod p$ with $i<p-m$ are thus obtained by scaling the columns of $G_{p,m}$ and deleting rows. They have rank at most $\min(16m,p)$ since $\rank_{\F_p}G_{p,m}\le p$. There are eight rows for each $i$ with $p-m\le i<2m$, hence $8(3m-p)_+$ remaining rows. Thus
\[
 \rho\le\min(16m,p)+8(3m-p)_+,\qquad
 v_p(\det K_m^{(p)})\ge16m-\rho\ge(16m-p)_+-8(3m-p)_+.
\]
The primes $p\le m$ contribute a nonnegative sum. Summing the bound for $p>m$ and applying \cref{lem:weighted-prime-sum} with $(a,b)=(1,16)$ and $(1,3)$ gives
\[
 \dfrac{1}{m^2}\sum_p v_p(\det K_m^{(p)})\log p
 \ge\dfrac{(16-1)^2-8(3-1)^2}{2}+o(1)
 =\dfrac{193}{2}+o(1).\qedhere
\]
\end{proof}

\begin{proposition}\label{prop:arithmetic}
Suppose that~\eqref{eq:vanish} holds for all sufficiently large primes. As $m=4\cdot17^r\to\infty$,
\[
 \log|\det K_m|\ge\left(-\frac{133}{4}+o(1)\right)m^2.
\]
\end{proposition}

\begin{proof}
By \cref{prop:matrix}, $\det K_m\in\Q^\times$. Since $\det K_m^{(p)}=p^{E_p}\det K_m$, the product formula and \cref{lem:total-exponents,lem:scaled-divisibility} give
\begin{align*}
 \frac{\log|\det K_m|}{m^2}
 &=\frac1{m^2}\sum_p\bigl(v_p(\det K_m^{(p)})-E_p\bigr)\log p\\*
 &\ge\frac{193}{2}-\frac{519}{4}+o(1)=-\frac{133}{4}+o(1).\qedhere
\end{align*}
\end{proof}

\section{Proofs of the main results}\label{sec:proofs}

\begin{proof}[Proof of \cref{thm:main}]
Suppose that $Z_\calA(3)=Z_\calA(5)=0$. Then~\eqref{eq:vanish} holds for all sufficiently large primes. Taking $m=4\cdot17^r$ in \cref{prop:matrix,prop:arithmetic} gives the contradiction
\begin{equation*}
 -\frac{133}{4}
 \le \liminf_{r\to\infty}\frac{\log|\det K_m|}{m^2}
 \le \limsup_{r\to\infty}\frac{\log|\det K_m|}{m^2}
 \le-24\log4<-\frac{133}{4}.\qedhere
\end{equation*}
\end{proof}

\begin{proof}[Proof of \cref{cor:large-weight}]
Fix an integer $N\ge3$. By \cref{thm:main}, there are infinitely many primes $p>2N+2$ for which at least one of $B_{p-3}$ and $B_{p-5}$ is non-zero modulo $p$. For each such prime, define in $\F_p[t]/(t^{2N+1})$
\[
 F_p(t)\coloneqq1+\sum_{n=1}^{p-1}\frac{t^2}{n^2}
    \prod_{j=1}^{n-1}(1+t/j)^2\prod_{j=n+1}^{p-1}(1-t/j)^2.
\]
The generating function of Maesaka and Ono~\cite{MO}*{Lemma~2.1 and Section~3, (1)} with $a=b=c=d=t$ gives
\[
 \log F_p(t)=\sum_{\substack{3\le h\le2N\\h\ \mathrm{odd}}}
    \frac{(2^h-2)^2}{h^2}B_{p-h}t^h.
\]
Here $\prod_{j=1}^{p-1}(1-t/j)=1-t^{p-1}\equiv1\pmod{t^{2N+1}}$ in $\F_p[t]$, and $B_{p-h}=0$ for even $2\le h\le2N$. Since $\log F_p(t)$ contains only odd powers of degree at least three, we have
\[
 F_p(t)F_p(-t)\equiv1\pmod{t^{2N+1}},\qquad
 [t^3]F_p(t)=4B_{p-3},\qquad [t^5]F_p(t)=36B_{p-5}.
\]
In the formula for $F_p(t)$, the factor $t^2/n^2$ contributes the exponent $2$ at $n$, and a nonconstant term of $(1\pm t/j)^2$ contributes an exponent $1$ or $2$ at $j$. Summing over the selected integers in increasing order expresses $[t^h]F_p(t)$ for $1\le h\le2N$ as a linear combination of the $\zeta_p(\bk)$ with entries in $\{1,2\}$ and weight $h$.

Put $P(t)\coloneqq\sum_{h=0}^{N}([t^h]F_p(t))t^h$. Its constant term is $P(0)=F_p(0)=1$. Suppose $\zeta_p(\bk)\equiv0\pmod p$ for every $\bk$ with entries in $\{1,2\}$ and $N<\wt(\bk)\le2N$. Then $[t^h]F_p(t)=0$ for $N<h\le2N$, so $P(t)\equiv F_p(t)\pmod{t^{2N+1}}$ and
\[
 P(t)P(-t)\equiv F_p(t)F_p(-t)\equiv1\pmod{t^{2N+1}}.
\]
Since $\deg(P(t)P(-t))=2\deg P\le2N$, we have $P(t)P(-t)=1$ in $\F_p[t]$. Hence $\deg P=0$ and $P(t)=1$. As $5\le2N$, both $[t^3]F_p(t)$ and $[t^5]F_p(t)$ would vanish, contrary to the choice of $p$.

Thus for each of the infinitely many primes $p$ above, some $\bk$ with entries in $\{1,2\}$ and $N<\wt(\bk)\le2N$ satisfies $\zeta_p(\bk)\not\equiv0\pmod p$. Only finitely many indices meet these conditions on the entries and weight, so one fixed $\bk$ satisfies the nonvanishing congruence for infinitely many $p$. Hence $\zeta_\calA(\bk)\ne0$.
\end{proof}

\subsection*{Acknowledgements}

The author was supported by JSPS KAKENHI (JP24K16901) and by access to OpenAI models provided through the ChatGPT for Academic Researchers program.

\subsection*{AI disclosure}

The author used ChatGPT to help develop the proofs, carry out computations, and prepare the manuscript. The author takes full responsibility for the content.

\end{document}